\input ./style/arxiv-general.cfg
\documentclass[MSNbibl,number,citesort,dvips]{arxbj}
\makeatletter
   \@ifpackageloaded{graphicx}{}{\usepackage{graphicx}}
\makeatother

\aid{0}
\volume{21}
\issue{3}
\pubyear{2015}
\firstpage{1361}
\lastpage{1385}
\doi{10.3150/14-BEJ605} 
\docsubty{FLA}

\makeatletter

\newtheorem{theorem}{Theorem}[section]
\newtheorem{proposition}[theorem]{Proposition}
\newtheorem{lemma}[theorem]{Lemma}
\newtheorem{corollary}[theorem]{Corollary}
\newremark{remark}[theorem]{Remark}

\newcommand{\cB}{\mathcal{B}}
\newcommand{\cI}{\mathcal{I}}
\newcommand{\cN}{\mathcal{N}}
\newcommand{\cX}{\mathcal{X}}

\newcommand{\Real}{\mathbb{R}}

\newcommand{\Esp}{\mathbb{E}}
\renewcommand{\Pr}{\mathbb{P}}

\renewcommand{\phi}{\varphi}
\renewcommand{\epsilon}{\varepsilon}
\renewcommand{\hat}{\widehat}
\newcommand{\eps}{\varepsilon}

\newcommand{\eqref}[1]{(\ref{#1})}

\def\vfrac#1#2{(#1)/#2}

\def\sklfrac#1#2{(#1/#2)}
\def\sklvfrac#1#2{((#1)/#2)}
\def\sklafrac#1#2{(#1/(#2))}

\makeatother

\begin{document}
\begin{frontmatter}

\title{Concentration inequalities for sampling
without replacement}
\runtitle{Concentration inequalities for sampling
without replacement}

\begin{aug}
\author[1]{\inits{R.}\fnms{R\'emi} \snm{Bardenet}\corref{}\thanksref{1}\ead[label=e1]{remi.bardenet@gmail.com}}
\and
\author[2]{\inits{O.-A.}\fnms{Odalric-Ambrym} \snm{Maillard}\thanksref{2}\ead[label=e2]{odalric.maillard@ee.technion.ac.il}}
\address[1]{Department of Statistics, University of Oxford, 1 South
Parks Road, OX1 3TG Oxford, UK.\\ \printead{e1}}
\address[2]{Faculty of Electrical Engineering, The Technion, Fishbach
Building, 32000 Haifa, Israel.\\ \printead{e2}}
\end{aug}

\received{\smonth{9} \syear{2013}}
\revised{\smonth{1} \syear{2014}}

%
\begin{abstract}
Concentration inequalities quantify the deviation of a random variable
from a fixed value. In spite of numerous applications, such as opinion
surveys or ecological counting procedures, few concentration results
are known for the setting of sampling without replacement from a
finite population. Until now, the best general concentration
inequality has been a Hoeffding inequality due to Serfling
[\textit{Ann. Statist.} \textbf{2} (1974) 39--48]. In
this paper, we first improve on the fundamental result
of Serfling [\textit{Ann. Statist.} \textbf{2} (1974)
39--48], and further extend it to
obtain a Bernstein
concentration bound for sampling without
replacement. We then derive an empirical version of our bound that does
not require the variance to be known to the user.
\end{abstract}

%
\begin{keyword}
\kwd{Bernstein}
\kwd{concentration bounds}
\kwd{sampling without replacement}
\kwd{Serfling}
\end{keyword}
\end{frontmatter}

\section{Introduction}
Few results exist on the concentration properties of sampling without
replacement from a finite population $\cX$. However, potential
applications are numerous, from historical applications such as
opinion surveys (Kish \cite{Kis65}) and ecological
counting procedures (Bailey \cite{Bai51}), to more recent
approximate Monte
Carlo Markov chain algorithms that use subsampled likelihoods
(Bardenet, Doucet and Holmes \cite{BaDoHo14}).
In a fundamental paper on sampling without
replacement, Serfling \cite{Serfling} introduced an
efficient Hoeffding bound,
that is, one which is a function of the range of the
population. Bernstein bounds are typically tighter when the variance
of the random variable under consideration is small, as
their leading term is linear in the standard deviation of
$\cX$, while the range only influences higher-order terms. This paper
is devoted to Hoeffding and Bernstein bounds for sampling without
replacement.


\subsection*{Setting and notations}
Let $\cX=(x_1,\dots,x_N)$ be a finite
population of $N$ real points. We use capital letters to denote random
variables on $\cX$, and lower-case letters for their possible
values. Sampling without replacement a list $(X_1,\dots,X_n)$ of size
$n$ from $\cX$ can be described sequentially as follows:
let first $\cI_1=\{1,\dots,n\}$, sample an integer $I_1$ uniformly on
$\cI_1$,
and set $X_1$ to be $x_{I_1}$. Then, for
each $i=2,\dots,n$, sample $I_i$ uniformly on the remaining indices
$\cI_{i} = \cI_{i-1}\setminus\{I_{i-1}\}$. Hereafter, we assume that
$N\geq2$.

\subsection*{Previous work}
There have been a few papers on concentration properties of sampling
without replacement; see, for instance, Hoeffding \cite{Hoeffding}, Serfling \cite{Serfling}, Horvitz and Thompson \cite{Horvitz52}, McDiarmid \cite{McDiarmid97}.
One notable contribution is the following reduction result in
Hoeffding's seminal paper
(Hoeffding \cite{Hoeffding}, Theorem~4):
%
\begin{lemma}\label{lem:reduction}
Let $\cX=(x_1,\ldots,x_N)$ be a finite population of $N$ real points,
$X_1,\ldots,X_n$ denote a random sample without replacement from $\cX$ and
$Y_1,\ldots,Y_n$ denote a random sample with replacement from $\cX$. If
$f\dvtx \Real\to\Real$ is continuous and convex, then
\begin{eqnarray*}
\Esp f \Biggl(\sum_{i=1}^n
X_i \Biggr) \leq\Esp f \Biggl(\sum_{i=1}^n
Y_i \Biggr) .
\end{eqnarray*}
\end{lemma}

Lemma~\ref{lem:reduction} implies that the concentration
results known for sampling with replacement as Chernoff bounds
(Boucheron, Lugosi and Massart \cite{BoLuMa13}) can be transferred to the
case of sampling without replacement. In particular, Proposition~\ref
{prop:hoeffnaive}, due to Hoeffding \cite{Hoeffding},
holds for the setting
without replacement.
%
\begin{proposition}[(Hoeffding's inequality)]\label{prop:hoeffnaive}
Let $\cX=(x_1,\ldots,x_N)$ be a finite population of $N$ points and
$X_1,\ldots,X_n$ be a random sample drawn without replacement from $\cX$.
Let
\[
a= \min_{1\leq i\leq N} x_i \quad \mbox{and}\quad b = \max
_{1\leq i\leq N} x_i.
\]
Then, for all $\epsilon>0$,
%
\begin{eqnarray}
\label{eqn:hoeffnaive} \Pr \Biggl(\frac{1}{n}\sum_{i=1}^n
X_i -\mu\geq\epsilon \Biggr) \leq\exp \biggl( - \frac{2n\epsilon^2}{(b-a)^2}
\biggr) ,
\end{eqnarray}
where $\mu= \frac{1}{N}\sum_{i=1}^N x_i$ is the mean of $\cX$.
\end{proposition}

The proof of Proposition~\ref{prop:hoeffnaive} (see, e.g.,
Boucheron, Lugosi and Massart \cite{BoLuMa13}) relies on a classical
bound on the moment-generating function of a random variable, which
we restate here as Lemma~\ref{l:hoeffnaive} for further reference.
%
\begin{lemma}
\label{l:hoeffnaive}
Let $X$ be a real random variable such that $\mathbb{E}X=0$ and
$a\leq X\leq b$ for some $a,b\in\mathbb{R}$. Then, for all $s \in
\mathbb{R}$,
\[
\log\mathbb{E}\mathrm{e}^{sX}\leq\frac{s^2(b-a)^2}{8}.
\]
\end{lemma}

When the variance of $\cX$ is small compared to the range $b-a$,
another Chernoff bound, known as Bernstein's bound (Boucheron, Lugosi and Massart \cite{BoLuMa13}),
is usually
tighter than Proposition~\ref{prop:hoeffnaive}.
%
\begin{proposition}[(Bernstein's inequality)]\label{prop:bernnaive}
With the notations of Proposition~\ref{prop:hoeffnaive}, let
\[
\sigma^2 = \frac{1}{N}\sum_{i=1}^N
(x_i-\mu)^2
\]
be
the variance of $\cX$. Then, for all $\epsilon>0$,
\begin{eqnarray*}
\Pr \Biggl(\frac{1}{n}\sum_{i=1}^n
X_i -\mu\geq\epsilon \Biggr) \leq\exp \biggl( - \frac{n \epsilon^2 }{2\sigma^2 + \sklfrac
{2}{3}(b-a)\epsilon}
\biggr) .
\end{eqnarray*}
\end{proposition}

Although these are interesting results, it appears that the bounds in
Propositions \ref{prop:hoeffnaive} and
\ref{prop:bernnaive} are actually very conservative,
especially when $n$ is large, say,
$n\geq N/2$. Indeed, Serfling \cite{Serfling} proved
that the term $n$ in the
RHS of \eqref{eqn:hoeffnaive} can be replaced
by $\frac{n}{1- (n-1)/N}$; see Theorem~\ref{thm:Hoeff}\vspace*{2pt} below, where the
result of Serfling is restated in our notation and slightly
improved. As $n$ approaches $N$, the bound of Serfling
\cite{Serfling}
improves dramatically, which corresponds to the intuition that when
sampling without replacement, the sample mean becomes a very accurate
estimate of $\mu$ as $n$ approaches $N$.

\subsection*{Contributions and outline} In Section~\ref{sec:Hoeff}, we
slightly modify Serfling's result, yielding a Hoeffding--Serfling bound in
Theorem~\ref{thm:Hoeff} that dramatically
improves on Hoeffding's in Proposition~\ref{prop:hoeffnaive}. In
Section~\ref{sec:Bern}, we contribute in Theorem~\ref{thm:Bern2} a similar
improvement on Proposition~\ref{prop:bernnaive}, which we call a
Bernstein--Serfling bound. To allow practical applications of our
Bernstein--Serfling bound, we finally provide an \emph{empirical}
Bernstein--Serfling bound in Section~\ref{sec:empBern}, in the spirit
of Maurer and Pontil \cite{MaurerP09}, which does not require
the variance of $\cX$ to
be known beforehand. In Section~\ref{sec:discussion}, we discuss
direct applications and potential further improvements of our results.

%
\begin{figure}

\includegraphics{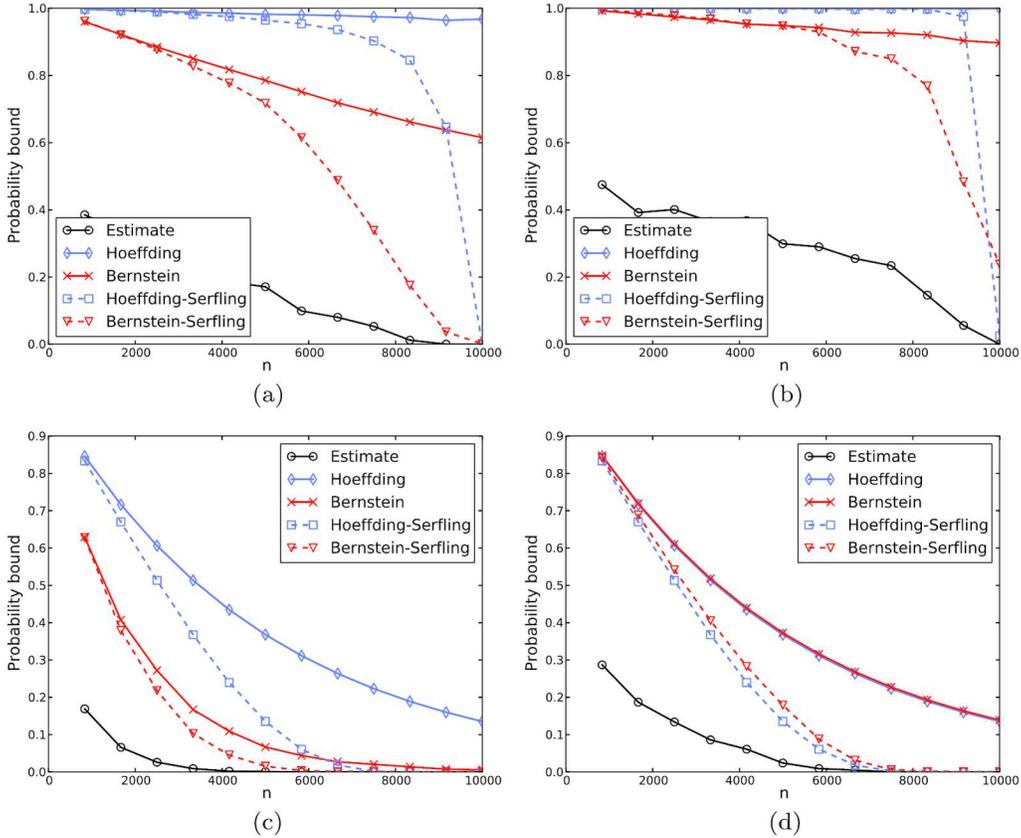}

\caption{Comparing known bounds on $p=\Pr(n^{-1}\sum_{i=1}^n X_i
-\mu
\geq0.01)$ with
our Hoeffding--Serfling and Bernstein--Serfling bounds. $\cX$ is here
a sample of size $N=10^4$ from each of the four distributions written
below each plot. An estimate (black plain line) of $p$
is obtained by averaging over $1000$ repeated subsamples of size $n$,
taken from $\cX$ uniformly without replacement. (a) Gaussian $\cN
(0,1)$. (b)
Log-normal $\ln\cN(1,1)$. (c) Bernoulli $\cB(0.1)$. (d) Bernoulli
$\cB(0.5)$.}\label{fig:plots}
\end{figure}

\subsection*{Illustration}
To give the reader a visual intuition of how the above mentioned bounds
compare in
practice and motivate their derivation, in Figure~\ref{fig:plots}, we
plot the bounds given by Proposition~\ref{prop:hoeffnaive} and
Theorem~\ref{thm:Hoeff} for Hoeffding bounds, and
Proposition~\ref{prop:bernnaive} and Theorem~\ref{thm:Bern2} for
Bernstein bounds for $\eps=10^{-2}$, in some common situations. We set
$\cX$ to be an
independent sample of size $N=10^4$ from each of the following four
distributions: unit centered Gaussian, log-normal with parameters
$(1,1)$, and Bernoulli with parameter 1$/$10 and 1$/$2. An estimate of the
probability
$\Pr(n^{-1}\sum_{i=1}^n X_i -\mu\geq10^{-2})$ is
obtained by averaging over $1000$ repeated samples of size $n$ taken
without replacement. In Figures~\ref{fig:plots}(a),
\ref{fig:plots}(b) and \ref{fig:plots}(c), Hoeffding's
bound and the Hoeffding--Serfling bound of Theorem~\ref{thm:Hoeff} are
close for $n\leq N/2$, after which the Hoeffding--Serfling bound decreases
to zero, outperforming Hoeffding's bound. Bernstein's and
our Bernstein--Serfling bound behave similarly, both outperforming their
counterparts that do not make use of the variance of
$\cX$. However, Figure~\ref{fig:plots}(d) shows that one
should not always prefer Bernstein bounds. In this case, the
standard deviation is as large as roughly half the range, making
Hoeffding's and Bernstein's bounds identical, and Hoeffding--Serfling
actually slightly better than Bernstein--Serfling. We emphasize here
that Bernstein bounds are typically useful when the variance is small
compared to the range.

\section{A reminder of Serfling's fundamental result}\label{sec:Hoeff}
In this section, we recall an initial result and proof by
Serfling \cite{Serfling}, and slightly improve on his
final bound.

We start by identifying the following martingales structures.
Let us introduce, for $1\leq k\leq N$,
%
\begin{eqnarray}\label{def:ZZs}
Z_k = \frac{1}{k}\sum_{t=1}^k
(X_t - \mu)\quad  \mbox{and}\quad  Z_k^\star=
\frac{1}{N-k}\sum_{t=1}^k
(X_t - \mu),\qquad  \mbox{where } \mu= \frac
{1}{N}\sum
_{i=1}^Nx_i .
\end{eqnarray}
Note that by definition $Z_N=0$, so that the $\sigma$-algebra
$\sigma(Z_{k+1},\ldots,Z_N)$ is equal to $\sigma(Z_{k+1},\ldots,\allowbreak Z_{N-1})$.
%
\begin{lemma}\label{lem:martingales}
The following \textit{forward martingale} structure holds for $\{
Z_k^\star\}_{k\leq N}$:
%
\begin{eqnarray}\label{eqn:fwmart}
\Esp \bigl[Z_k^\star | Z_{k-1}^\star,
\ldots,Z_{1}^\star \bigr] = Z_{k-1}^\star.
\end{eqnarray}
Similarly, the following \textit{reverse martingale} structure holds
for $\{Z_k\}_{k\leq N}$:
%
\begin{eqnarray}\label
{eqn:rvmart}
\Esp [Z_k | Z_{k+1},\ldots,Z_{N-1} ] =
Z_{k+1} .
\end{eqnarray}
\end{lemma}
%
%
\begin{pf}
We first prove \eqref{eqn:fwmart}. Let $1\leq k \leq N$. We start by
noting that
%
\begin{eqnarray}\label
{eqn:decZ*}
Z_k^\star&=& \frac{1}{N-k}\sum
_{t=1}^{k-1} (X_t - \mu) +
\frac
{X_k - \mu}{N-k}
\nonumber
\\[-8pt]\\[-8pt]
&=& \frac{N-k+1}{N-k}Z_{k-1}^\star+ \frac{X_k - \mu}{N-k}
.\nonumber
\end{eqnarray}
%
Since $X_k$ is uniformly distributed on the remaining elements of
$\cX$ after $X_1,\ldots,X_{k-1}$ have been drawn, its conditional
expectation given $X_1,\ldots,X_{k-1}$ is the average of the $N-k+1$
remaining points in $\cX$. Since points in $\cX$ add up to $N\mu$,
we obtain
%
\begin{eqnarray}\label{eqn:EspXk}
\Esp \bigl[ X_k |Z_{k-1}^\star,
\ldots,Z_{1}^\star \bigr] &=& \Esp [ X_k
|X_{k-1},\ldots,X_{1} ]
\nonumber
\\
&=& \frac{N\mu- \sum_{i=1}^{k-1}X_i}{N-k+1}
\\
&=& \mu- Z_{k-1}^\star.\nonumber
\end{eqnarray}
Combined with \eqref{eqn:decZ*}, this yields \eqref{eqn:fwmart}.

We now turn to proving \eqref{eqn:rvmart}. First, let $1\leq k\leq
N$. Since
\[
(k+1)Z_{k+1}= (N-k-1)\mu-X_{k+2}-\cdots-X_N,
\]
$\sigma(Z_{k+1},\ldots,Z_{N-1})$ is equal to
$\sigma(X_{k+2},\ldots,X_N)$. Now, let us remark that
the indices of $(X_1,\ldots,X_N)$ are uniformly distributed on the
permutations of
$\{1,\ldots,N\}$, so that $(X_1,\ldots,\allowbreak X_{N-k})$ and
$(X_{k+1},\ldots,X_N)$ have the same marginal distribution. Consequently,
\[
\Esp [X_{k+1} \vert Z_{k+1},\ldots,Z_{N-1} ] = \Esp
[X_{k+1} \vert X_{k+2}\ldots, X_{N} ] =
\frac{S_{k+1}}{k+1},
\]
where we introduced the sum $S_{k+1} = \sum_{t=1}^{k+1}X_t$.
Finally, we prove \eqref{eqn:rvmart} along the same lines as
\eqref{eqn:fwmart}:
\begin{eqnarray*}
\Esp [Z_k \vert Z_{k+1},\ldots,Z_{N-1} ] &=& \Esp
\biggl[\frac{S_k-k\mu}{k} \Bigl\vert Z_{k+1},\ldots,Z_{N-1} \biggr]
\\
&=& \Esp \biggl[\frac{S_{k+1}-X_{k+1}}{k} \Bigl\vert Z_{k+1},\ldots,Z_{N}
\biggr] -\mu
\nonumber
\\
&=&\frac{S_{k+1}}{k} - \frac{S_{k+1}}{k(k+1)} -\mu
\nonumber
\\
&=& Z_{k+1}.
\nonumber
\end{eqnarray*}
\upqed
\end{pf}

\subsection*{A Hoeffding--Serfling inequality}
Let us now state the main result of Serfling \cite{Serfling}. This is a key
result to derive a concentration inequality, a maximal
concentration inequality and a self-normalized concentration
inequality, as explained in Serfling \cite{Serfling}.
%
\begin{proposition}[(Serfling \cite{Serfling})]\label
{prop:Hoeff}
Let us denote $a= \min_{1\leq i\leq N} x_i$, and $b = \max_{1\leq
i\leq N} x_i$.
Then, for any $\lambda>0$, it holds that
\begin{eqnarray*}
\log\Esp\exp ( \lambda n Z_n ) \leq\frac{(b-a)^2}{8}
\lambda^2 n \biggl(1- \frac{n-1}{N} \biggr) .
\end{eqnarray*}
Moreover, for any $\lambda>0$, it also holds that
\begin{eqnarray*}
\log\Esp\exp \Bigl( \lambda\max_{1\leq k\leq n} Z_k
^\star \Bigr) \leq\frac{(b-a)^2}{8} \frac{\lambda^2}{(N -n)^2} n \biggl(1-
\frac
{n-1}{N} \biggr) .
\end{eqnarray*}
\end{proposition}
\begin{pf}
First, \eqref{eqn:decZ*} yields that for all $\lambda'>0$,
%
\begin{eqnarray}\label{eqn:th1dec}
\lambda' Z_{k}^\star= \lambda'
Z_{k-1}^\star+\lambda' \frac
{X_k-\mu+ Z_{k-1}^\star}{N-k}
.
\end{eqnarray}
Furthermore, we know from \eqref{eqn:EspXk} that
$-Z_{k-1}^\star$ is
the conditional expectation of $X_k-\mu$ given
$X_1,\ldots,X_{k-1}$. Thus, since $X_k- \mu\in
[a-\mu,b-\mu]$, Lemma~\ref{l:hoeffnaive} applies and we get
that, for all $2\leq k\leq n$,
%
\begin{eqnarray}\label{eqn:th1Hoef}
\log\Esp \biggl[ \exp \biggl( \lambda' \frac{X_k-\mu+ Z_{k-1}^\star
}{N-k} \biggr) \Bigl|
Z_1^\star,\ldots,Z_{k-1}^\star \biggr]
\leq\frac{(b-a)^2}{8} \frac{{\lambda'}^2}{ (N-k )^2} .
\end{eqnarray}
Similarly, we can apply Lemma~\ref{l:hoeffnaive} to
$Z_1^\star= (X_1-\mu)/(N-1)$ to obtain
%
\begin{eqnarray}\label{eqn:th1Hoef2}
\log\Esp\exp \bigl( \lambda' Z_1^\star \bigr)
\leq\frac
{(b-a)^2}{8} \frac{{\lambda'}^2}{ (N-1 )^2} .
\end{eqnarray}
Upon noting that $Z_n = \frac{N-n}{n}Z_n^\star$, and combining \eqref
{eqn:th1Hoef} and \eqref{eqn:th1Hoef2} together
with the decomposition \eqref{eqn:th1dec}, we eventually obtain the bound
\begin{eqnarray*}
\log\Esp\exp \biggl( \lambda' \frac{n}{N-n} Z_n
\biggr) \leq\frac
{(b-a)^2}{8} \sum_{k=1}^n
\frac{{\lambda'}^2}{(N-k)^2} .
\end{eqnarray*}
In particular, for $\lambda$ such that $\lambda' = (N-n)\lambda$, the
RHS of this equation contains the quantity
%
\begin{eqnarray}\label{eqn:sumkn}
\sum_{k=1}^n \frac{(N-n)^2}{(N-k)^2} &=& 1 +
(N-n)^2\sum_{k=N-n+1}^{N-1}
\frac{1}{k^2}
\nonumber
\\
&\leq& 1+ (N-n)^2\frac{((N-1) - (N-n))}{(N-n)N} = 1+ \frac
{(N-n)(n-1)}{N}
\\
&=& 1 + n-1 - n \frac{n-1}{N} = n \biggl(1 - \frac{n-1}{N} \biggr)
,\nonumber
\end{eqnarray}
where we used in the second line the following approximation from
(Serfling \cite{Serfling}, Lemma~2.1): for
$1\leq j\leq m$, it holds
\begin{eqnarray*}
\sum_{k=j+1}^{l} \frac{1}{k^2} \leq
\frac{l-j}{j(l+1)} .
\end{eqnarray*}
This concludes the proof of the first result of
Proposition~\ref{prop:Hoeff}. The second result follows from applying
Doob's maximal inequality for martingales combined with the previous derivation.
\end{pf}

The result of Proposition~\ref{prop:Hoeff} reveals a powerful feature
of the no replacement setting: the factor
$n(1- \frac{n-1}{N})$ in the exponent, as opposed to $n$ in the case
of sampling with replacement. This leads to a dramatic improvement of
the bound when $n$ is large, as can be seen on
Figure~\ref{fig:plots}. Serfling \cite{Serfling}
mentioned that
a factor $1 - \frac{n}{N}$ would be intuitively more natural, as
indeed when $n=N$ the mean $\mu$ is known exactly, so that $Z_N$ is
deterministically zero.

Serfling did not publish any result with $1 - \frac{n}{N}$. However,
it appears that a careful examination of the previous proof and of the
use of equation \eqref{eqn:rvmart}, in lieu of
\eqref{eqn:fwmart}, allows us to get such an improvement. We detail
this in the following proposition. More than a simple cosmetic
modification, it is actually a slight improvement on Serfling's original
result when $n>N/2$.
%
\begin{proposition}\label{prop:Hoeff2}
Let $(Z_k)$ be defined by \eqref{def:ZZs}. For any $\lambda>0$, it
holds that
\begin{eqnarray*}
\log\Esp\exp ( \lambda n Z_n ) \leq\frac{(b-a)^2}{8}
\lambda^2(n+1) \biggl(1 - \frac{n}{N} \biggr) .
\end{eqnarray*}
Moreover, for any $\lambda>0$, it also holds that
\begin{eqnarray*}
\log\Esp\exp \Bigl( \lambda\max
_{n\leq k\leq N-1} Z_k \Bigr) \leq \frac{(b-a)^2}{8}
\frac{\lambda^2}{n^2}(n+1) \biggl(1 - \frac
{n}{N} \biggr) .
\end{eqnarray*}
\end{proposition}
\begin{pf}
Let us introduce the notation $Y_k = Z_{N-k}$ for $1\leq k\leq N-1$.
From \eqref{eqn:rvmart}, it comes
\begin{eqnarray*}
\Esp [Y_{N-k} | Y_{1},\ldots, Y_{N-k-1} ] =
Y_{N-k-1} .
\end{eqnarray*}
By a change of variables,
this can be rewritten as
\begin{eqnarray*}
\Esp [Y_{k} | Y_{1},\ldots, Y_{k-1} ] =
Y_{k-1} .
\end{eqnarray*}

Now we remark that the following decomposition holds:
%
\begin{eqnarray}\label{eqn:Ydec}
\lambda Y_k &=& \lambda\frac{\sum_{i=1}^{N-k}(X_i-\mu
)}{N-k}
\nonumber
\\[-8pt]\\[-8pt]
&=& \lambda Y_{k-1} - \lambda\frac{X_{N-k+1}-\mu-Y_{k-1}}{N-k} .\nonumber
\end{eqnarray}
Since $Y_{k-1}$ is the conditional mean of
$X_{N-k+1}-\mu\in[a-\mu,b-\mu]$, Lemma~\ref{l:hoeffnaive} yields that,
for all $2\leq k\leq n$,
%
\begin{eqnarray}\label{eqn:th1Hoefb}
\log\Esp \biggl[ \exp \biggl( \lambda' \frac{X_{N-k+1}-\mu-
Y_{k-1}}{N-k} \biggr) \Bigl|
Y_1,\ldots,Y_{k-1} \biggr] \leq\frac{(b-a)^2}{8}
\frac{{\lambda'}^2}{ (N-k )^2} .
\end{eqnarray}
On the other hand it holds by definition of $Y_1$ that
\begin{eqnarray*}
Y_1 = Z_{N-1} = \frac{\sum_{i=1}^{N-1}(X_i-\mu)}{N-1} \in[a-\mu ,b-\mu] .
\end{eqnarray*}
Along the lines of the proof of Proposition~\ref{prop:Hoeff}, we obtain
%
\begin{eqnarray}\label{eqn:th1Hoef2b}
\log\Esp\exp \bigl( \lambda' Y_1 \bigr) \leq
\frac{(b-a)^2}{8} \frac{{\lambda'}^2}{ (N-1 )^2} .
\end{eqnarray}
Combining equations \eqref{eqn:th1Hoefb} and \eqref{eqn:th1Hoef2b}
with the decomposition \eqref{eqn:Ydec}, it comes
\begin{eqnarray*}
\log\Esp\exp \bigl( \lambda' Y_n \bigr) &\leq&
\frac{(b-a)^2}{8} \sum_{k=1}^n
\frac{{\lambda'}^2}{(N-k)^2}
\\
&\leq& \frac{(b-a)^2}{8} \frac{{\lambda'}^2}{(N-n)^2}n \biggl(1-\frac
{n-1}{N}
\biggr) ,
\end{eqnarray*}
where in the last line we made use of \eqref{eqn:sumkn}. Rewriting
this inequality in terms of $Z$, we obtain that, for all $1\leq n\leq N-1$,
\begin{eqnarray*}
\log\Esp\exp \bigl( \lambda(N-n) Z_{N-n} \bigr) \leq\frac
{(b-a)^2}{8}
\lambda^2n \biggl(1-\frac{n-1}{N} \biggr) ,
\end{eqnarray*}
that is, by resorting to a new change of variable,
\begin{eqnarray*}
\log\Esp\exp ( \lambda n Z_{n} ) &\leq& \frac{(b-a)^2}{8}
\lambda^2(N-n) \biggl(1-\frac{N-n-1}{N} \biggr)
\\
&\leq& \frac{(b-a)^2}{8} \lambda^2(N-n)\frac{n+1}{N}
\\
&\leq& \frac{(b-a)^2}{8} \lambda^2(n+1) \biggl(1 -
\frac{n}{N} \biggr) .
\end{eqnarray*}
The second part of the proposition follows from applying Doob's maximal
inequality for martingales to $Y_n$, similarly to Proposition~\ref{prop:Hoeff}.
\end{pf}
%
\begin{theorem}[(Hoeffding--Serfling inequality)]\label{thm:Hoeff}
Let $\cX=(x_1,\ldots,x_N)$ be a finite population of $N>1$ real points,
and $(X_1,\ldots,X_n)$ be a list of size $n<N$ sampled without
replacement from~$\cX$.
Then for all $\epsilon>0$, the following concentration bounds hold
\begin{eqnarray*}
\Pr \biggl( \max_{n\leq k\leq N-1} \frac{\sum_{t=1}^k (X_t - \mu
)}{k} \geq\eps \biggr) &
\leq& \exp \biggl( - \frac{2 n\epsilon
^2}{(1- n/N)(1+1/n)(b-a)^2} \biggr),
\\
\Pr \biggl( \max_{1\leq k\leq n} \frac{\sum_{t=1}^k (X_t - \mu
)}{N-k} \geq
\frac{n\eps}{N-n} \biggr) &\leq& \exp \biggl( - \frac
{2n\epsilon^2}{(1- (n-1)/N)(b-a)^2} \biggr) ,
\end{eqnarray*}
where $a= \min_{1\leq i\leq N} x_i$ and $b=\max_{1\leq i\leq N} x_i$.
\end{theorem}
\begin{pf}
Applying Proposition~\ref{prop:Hoeff2} together with Markov's
inequality, we obtain that, for all $\lambda>0$,
\begin{eqnarray*}
\Pr \biggl( \max_{n\leq k\leq N-1} \frac{\sum_{t=1}^k (X_t - \mu
)}{k} \geq\epsilon \biggr)
\leq \exp \biggl( -\lambda\epsilon+ \frac{(b-a)^2}{8} \frac{\lambda
^2}{n^2} (n+1)
(1- n/N) \biggr) .
\end{eqnarray*}
%
We now optimize the previous bound in $\lambda$. The optimal value is
given by
\begin{eqnarray*}
\lambda^\star= \epsilon\frac{4}{(b-a)^2} \frac{n^2}{(n+1)(1- n/N)} .
\end{eqnarray*}
This gives the first inequality of Theorem~\ref{thm:Hoeff}. The proof
of the second
inequality follows the very same lines.
\end{pf}

Inverting the result of Theorem~\ref{thm:Hoeff} for $n<N$ and
remarking that the
resulting bound still holds for $n=N$, we straightforwardly
obtain the following result.
%
\begin{corollary}\label{cor:Hoeff}
For all $n\leq N$, for all $\delta\in[0,1]$, with probability higher
than $1-\delta$, it holds
\begin{eqnarray*}
\frac{\sum_{t=1}^n (X_t - \mu)}{n} &\leq& (b-a)\sqrt{\frac{\rho
_n\log(1/\delta) }{2n}} ,
\end{eqnarray*}
where we define
%
\begin{eqnarray}\label{eqn:rhon}
\rho_n = %
\cases{ \biggl(1-\displaystyle \frac{n-1}{N}\biggr) &\quad
\mbox{if} $n\leq N/2$,
\cr
\biggl(1-\displaystyle \frac{n}{N}\biggr) (1+1/n) & \quad \mbox{if} $n
> N/2$ .} %
\end{eqnarray}
%
\end{corollary}

\section{A Bernstein--Serfling inequality}\label{sec:Bern}

In this section, we consider $\sigma^2 = N^{-1}\sum_{i=1}^N
(x_i-\mu)^2$ is known, and extend Theorem~\ref{thm:Hoeff} to that situation.

Similarly to Lemma~\ref{lem:martingales}, the following structural
lemma will be useful:
%
\begin{lemma}\label{lem:var}
It holds
\begin{eqnarray*}
\Esp \bigl[ (X_k - \mu)^2 | Z_1,\ldots,
Z_{k-1} \bigr] =\sigma^2 - Q_{k-1}^\star,
\qquad \mbox{where } Q_{k-1}^\star= \frac{\sum_{i=1}^{k-1}
 ((X_i-\mu)^2- \sigma^2 )}{N-k+1} ,
\end{eqnarray*}
where the $Z_i$'s are defined in \eqref{def:ZZs}. Similarly, it holds
\begin{eqnarray*}
\Esp \bigl[ (X_{k+1} - \mu)^2 | Z_{k+1},\ldots,
Z_{N-1} \bigr] =\sigma^2 + Q_{k+1},\qquad  \mbox{where }
Q_{k+1} = \frac{\sum_{i=1}^{k+1}
 ((X_i-\mu)^2- \sigma^2 )}{k+1} .
\end{eqnarray*}
\end{lemma}

\begin{pf}
We simply remark again that, conditionally on $X_1,\ldots,X_{k-1}$, the
variable $X_k$ is
distributed uniformly over the remaining points in $\cX$, so that
\begin{eqnarray*}
\Esp \bigl[ (X_k-\mu)^2 | Z_1,\ldots,
Z_{k-1} \bigr] &=& \Esp \bigl[ (X_k-\mu)^2 |
X_1,\ldots, X_{k-1} \bigr]
\\
&=& \frac{1}{N-k+1} \Biggl[ N\sigma^2 - \sum
_{i=1}^{k-1}(X_i-\mu )^2
\Biggr]
\\
&=& \sigma^2-Q_{k-1}^\star. 
\end{eqnarray*}
The second equality of Lemma~\ref{lem:var} follows from the same
argument, as in the proof of Lemma~\ref{lem:martingales}.
\end{pf}

Let us now introduce the following notations:
\begin{eqnarray*}
\mu_{<,k+1} &=& \Esp [X_{k+1} - \mu | Z_{k+1},\ldots,
Z_{N-1} ] ,
\\
\mu_{>,k} &=& \Esp [X_{k} - \mu | Z_{1},\ldots,
Z_{k-1} ] ,
\\
\sigma_{<,k+1}^2 &=& \Esp \bigl[ (X_{k+1} -
\mu)^2 | Z_{k+1},\ldots, Z_{N-1} \bigr] -
\mu_{<,k+1}^2 ,
\\
\sigma_{>,k}^2 &=& \Esp \bigl[ (X_{k} -
\mu)^2 | Z_{1},\ldots, Z_{k-1} \bigr] -
\mu_{>,k}^2 .
\end{eqnarray*}
%
We are now ready to state Proposition~\ref{prop:Bern}, which is a
Bernstein version of Proposition~\ref{prop:Hoeff}.

\begin{proposition}\label{prop:Bern}
For any $\lambda>0$, it holds that
\begin{eqnarray*}
\log\Esp\exp \Biggl( \lambda n Z_n - {\lambda}^2\sum
_{k=1}^{N-n} \phi \biggl(\frac{2(b-a) \lambda}{N-k}
\biggr)\frac{\sigma_{<,N-k+1}^2
n^2}{(N-k)^2} \Biggr) &\leq&0 ,
\\
\log\Esp\exp \Biggl( \lambda n Z_n - {\lambda}^2\sum
_{k=1}^{n} \phi \biggl(2(b-a) \lambda
\frac{N-n}{N-k} \biggr)\frac{\sigma_{>,k}^2
(N-n)^2}{(N-k)^2} \Biggr) &\leq&0 ,
\end{eqnarray*}
%
where we introduced the function $\phi(c) = \frac{\mathrm{e}^{c} -1 -c}{c^2}$.
Moreover, for any $\lambda>0$, it also holds that
\begin{eqnarray*}
\log\Esp\exp \Biggl( \lambda \Bigl(\max_{1\leq k\leq n}
Z_k^\star \Bigr) - \sum_{k=1}^{n}
\phi \biggl(\frac{2(b-a) \lambda}{N-k} \biggr)\frac{\sigma_{>,k}^2 {\lambda}^2}{(N-k)^2} \Biggr) &\leq&0 ,
\\
\log\Esp\exp \Biggl( \lambda \Bigl(\max_{n\leq k\leq N-1} Z_{k}
\Bigr) - \sum_{k=1}^{N-n} \phi \biggl(
\frac{2(b-a) \lambda}{N-k} \biggr)\frac
{\sigma_{<,N-k+1}^2 {\lambda}^2}{(N-k)^2} \Biggr) &\leq&0 .
\end{eqnarray*}
\end{proposition}
\begin{pf}
The key point is to replace equations \eqref{eqn:th1Hoef} and
\eqref{eqn:th1Hoef2} in the proof of Proposition~\ref{prop:Hoeff},
which make use of the range of $\cX$, by equivalent ones that involve the
variance. We only detail the proof of the first inequality, the proof
of the three others follows the same steps.

A standard result from the proof of Bennett's inequality (see
Lugosi \cite{Lug09}, page 11, or Boucheron, Lugosi and Massart
\cite{BoLuMa13}, proof of Theorem~2.9)
applied to the random variable
$X_{N-k+1}-\mu$, with conditional mean $\mu_{<,N-k+1}$ and
conditional variance $\sigma_{<,N-k+1}^2$, yields
%
\begin{eqnarray}\label{eqn:th2Ben}
&&\Esp \biggl[ \exp \biggl( \lambda' \frac{X_{N-k+1}- \mu+
Y_{k-1}}{N-k} \nonumber\\[-8pt]\\[-8pt]
&&\hphantom{\Esp \biggl[ \exp \biggl(}{} -
\sigma_{<,N-k+1}^2 \phi \biggl(\frac{2(b-a)\lambda'}{N - k} \biggr)
\frac{{\lambda'}^2}{ (N - k )^2} \biggr) \Bigl| Y_1,\ldots ,Y_{k-1} \biggr]
\leq1 ,\nonumber
\end{eqnarray}
where we used the notation $Y_{k}= Z_{N-k}$ of Proposition~\ref
{prop:Hoeff2}, and the function $\phi$ defined in the statement of the
proposition
\begin{eqnarray*}
\phi(c) = \frac{\mathrm{e}^{c} -1 -c}{c^2} .
\end{eqnarray*}
Similarly, $Y_1$ satisfies
%
\begin{eqnarray}\label{eqn:th2Ben2}
\log\Esp\exp \bigl( \lambda' Y_1 \bigr) = \log\Esp\exp
\biggl( \lambda' \frac{\mu-X_N}{N-1} \biggr) \leq
\sigma_{<,N}^2 \phi \biggl(\frac{2(b-a) \lambda'}{N-1} \biggr)
\frac{{\lambda'}^2}{ (N-1 )^2} ,
\end{eqnarray}
where $\sigma_{<,N}^2 = \sigma^2$ is deterministic.

Thus, combining \eqref{eqn:th2Ben} and \eqref{eqn:th2Ben2} together
with the decomposition \eqref{eqn:Ydec}, we eventually get the bound
\begin{eqnarray*}
\log\Esp\exp \Biggl( \lambda' Y_n- \sum
_{k=1}^n \phi \biggl(\frac
{2(b-a) \lambda'}{N-k} \biggr)
\frac{\sigma_{<,N-k+1}^2 {\lambda
'}^2}{(N-k)^2} \Biggr) &\leq& 0 . 
\end{eqnarray*}
\upqed
\end{pf}

Using the result of Proposition~\ref{prop:Bern}, we could immediately
derive a simple Bernstein inequality for
sampling without replacement via an application of Theorem~\ref{thm:Hoeff}
to the random variables $Z_i=(X_i-\mu)^2$.
However, Maurer and Pontil \cite{MaurerP09} and Audibert, Munos and Szepesv{\'a}ri \cite{AuMuSz09} showed that, in the case
of sampling with replacement, a careful use of self-bounded properties
of the
variance yields better bounds. We now explain how to get a similar
improvement on the naive Bernstein inequality in the case of
sampling without replacement. We start with a technical lemma.
%
\begin{lemma}\label{lem:sigmaleft}
For all $\delta\in[0,1]$, with probability larger than $1-\delta$,
it holds
%
\begin{eqnarray}\label{eqn:sigmaleft}
&&\max_{1\leq k\leq n}\sigma_{>,k}^2 \leq
\sigma^2 + \frac{\sigma(b-a)(n-1)}{N-n+1}\sqrt{\frac{2\log(1/\delta)}{n-1}} .
\end{eqnarray}
Similarly, with probability larger than $1-\delta$, it holds
%
\begin{eqnarray}\label{eqn:sigmaleft2}
&&\max_{n\leq k\leq N-1}\sigma_{<,k+1}^2 \leq
\sigma^2 + \frac
{\sigma(b-a)(N-n-1)}{n+1}\sqrt{\frac{2\log(1/\delta)}{N-n-1}} .
\end{eqnarray}
\end{lemma}
%
\begin{remark}
When $N\to\infty$, the upper bound on
$\max_{1\leq k\leq n}\sigma_{>,k}^2$ reduces to $\sigma^2$. Indeed,
this limit case intuitively corresponds to sampling with replacement,
for which the conditional variance equals $\sigma^2$.
\end{remark}
\begin{pf*}{Proof of Lemma~\ref{lem:sigmaleft}}
We first prove \eqref{eqn:sigmaleft}. By definition and Lemma~\ref
{lem:var}, it holds that
%
\begin{eqnarray}\label{eqn:decompsig2}
\sigma_{>,k}^2 &=& \sigma^2 -
Q_{k-1}^\star- {Z_{k-1}^\star}^2
\nonumber
\\[-8pt]\\[-8pt]
&\leq& \sigma^2 - \frac{1}{N-k+1}\sum
_{i=1}^{k-1} \bigl[ (X_i-\mu
)^2 - \sigma^2 \bigr] .\nonumber
\end{eqnarray}
Let $V_{k-1} = \frac{1}{k-1}\sum_{i=1}^{k-1}
 (X_i-\mu )^2$. Equation \eqref{eqn:decompsig2} yields
\begin{eqnarray*}
\max_{1\leq k\leq n} \sigma_{>,k}^2 \leq
\sigma^2 + \max_{1\leq
k\leq n}\frac{k-1}{N-k+1} \bigl(
\sigma^2 - V_{k-1} \bigr) .
\end{eqnarray*}
The rest of the proof proceeds by establishing a suitable
maximal concentration bound for the quantity $V_{k-1}$, the mean of
which is $\sigma^2$.

We remark that $-Q_{k-1}^\star= \frac{k-1}{N-k+1} (\sigma^2 -
V_{k-1} )$ is a martingale. Indeed, it satisfies
\begin{eqnarray*}
&&\Esp \bigl[ -Q^\star_{k-1} | Q^\star_{k-2},
\ldots,Q^\star _1 \bigr]
\\
&&\quad = \frac{1}{N-k+1}\Esp \Biggl[ \sum_{i=1}^{k-1}
\bigl(\sigma^2 - (X_i-\mu)^2 \bigr) |
Q^\star_{k-2},\ldots,Q^\star_1 \Biggr]
\\
&&\quad = \frac{1}{N-k+1}\sum_{i=1}^{k-2} \bigl(
\sigma^2 - (X_i-\mu)^2 \bigr) +
\frac{1}{N-k+1}\Esp \bigl[ \bigl(\sigma^2 - (X_{k-1}-
\mu)^2 \bigr) | Q^\star_{k-2},\ldots,Q^\star_1
\bigr]
\\
&&\quad = -\frac{N-k+2}{N-k+1}Q^\star_{k-2} + \frac{1}{N-k+1}Q_{k-2}^\star
\\
&&\quad = -Q^\star_{k-2} ,
\end{eqnarray*}
where we applied Lemma~\ref{lem:var} in the third line.
Doob's maximal inequality thus yields that, for all $\lambda>0$,
\begin{eqnarray*}
\Pr \Bigl(\max_{1\leq k\leq n} -Q^\star_{k-1} \geq
\epsilon \Bigr) &=& \Pr \Bigl(\max_{1\leq k\leq n} \exp\bigl(-\lambda
Q^\star_{k-1}\bigr) \geq\exp(\lambda\epsilon) \Bigr)
\\
&\leq& \Esp \bigl[\exp \bigl(-\lambda Q^\star_{n-1} - \lambda
\epsilon \bigr) \bigr]
\\
&=& \Esp \biggl[\exp \biggl(\lambda\frac{n-1}{N-n+1} \biggl(
\sigma^2 - V_{n-1}- \frac{N-n+1}{n-1} \epsilon \biggr)
\biggr) \biggr] .
\end{eqnarray*}

At this point, we fix $\lambda>0$ and apply Lemma~\ref{lem:reduction}
to the random variables
$X'_i = (X_i-\mu)^2$ and function $f\dvtx x\to\exp(-\lambda(n-1) x)$.
We deduce that, for all $\epsilon'>0$ and $\lambda>0$,\vspace*{-1pt}
%
\begin{eqnarray}\label{eqn:sigV1}
\Pr \biggl( \max_{1\leq k\leq n}\sigma_{>,k}^2 -
\sigma^2 \geq\frac
{n-1}{N-n+1}\epsilon' \biggr) &\leq&
\Esp \bigl[\exp \bigl(-\lambda\bigl(V_{n-1} - \sigma^2 +
\epsilon'\bigr) \bigr) \bigr]
\nonumber
\\[-8pt]\\[-8pt]
&\leq& \Esp \bigl[\exp \bigl(-\lambda\bigl(\tilde V_{n-1} -
\sigma^2 +\epsilon'\bigr) \bigr) \bigr] ,\nonumber
\end{eqnarray}
where we introduced in the last line the notation $\tilde V_{n-1} =
\frac{1}{n-1}\sum_{i=1}^{n-1}  (Y_i-\mu )^2$, with
the $\{Y_i\}_{1\leq i\leq n-1}$ being sampled from $\cX$ \textit{with
replacement}.
Note that $\tilde V_{n-1}$ has mean $\sigma^2$ too.

Now, we check that the assumptions of Theorem~13 of Maurer \cite{Maurer2006}
hold. We first introduce the modification\vspace*{-1pt}
\[
\mathbf{Y}^{j,y}_{1:n-1} = \{Y_1,
\ldots,Y_{j-1},y,Y_{j+1},\ldots, Y_{n-1}\}
\]
of $\mathbf{Y}_{1:n-1}$, where $Y_j$ is replaced by $y\in\cX$.
Writing $\tilde V_{n-1} =\tilde V_{n-1}(\mathbf{Y}_{1:n-1})$ to
underline the dependency on the sample set $\mathbf{Y}_{1:n-1}$,
it straightforwardly comes, on the one hand, that for all $y\in\cX$\vspace*{-1pt}
\begin{eqnarray*}
\tilde V_{n-1}(\mathbf{Y}_{1:n-1}) - \tilde V_{n-1}
\bigl(\mathbf {Y}^{j,y}_{1:n-1}\bigr) &=& \frac{1}{n-1}
\bigl((Y_j-\mu)^2 - (y-\mu )^2 \bigr)
\\
&\leq& \frac{1}{n-1} (Y_j-\mu)^2 \leq
\frac{1}{n-1}(b-a)^2 ,
\end{eqnarray*}
and, on the other hand, that the following self-bounded property holds:
\begin{eqnarray*}
\sum_{j=1}^{n-1} \Bigl(\tilde
V_{n-1}(\mathbf{Y}_{1:n-1}) - \inf_{y\in\cX}\tilde
V_{n-1}\bigl(\mathbf{Y}^{j,y}_{1:n-1}\bigr)
\Bigr)^2 &\leq &\frac{1}{(n-1)^2} \sum_{j=1}^{n-1}(Y_j-
\mu)^4
\\
&\leq& \frac{(b-a)^2}{n-1}\tilde V_{n-1}(\mathbf{Y}_{1:n-1}) .
\end{eqnarray*}
We now apply of the proof of
Theorem~13 of Maurer \cite{Maurer2006}\footnote{The
theorem is stated for
$\Pr [\Esp[Z] - Z \geq\epsilon ]$ but, actually,
$\Esp [\exp (-\lambda(Z - \Esp[Z] +\epsilon) ) ]$ is
bounded in the proof.} to $Z= \frac{n-1}{(b-a)^2}\tilde V_{n-1}$,
together with \eqref{eqn:sigV1}, which yields
\begin{eqnarray*}
\Pr \biggl( \max_{1\leq k\leq n}\sigma_{>,k}^2 -
\sigma^2 \geq\frac
{(b-a)^2}{N-n+1}\epsilon \biggr) &\leq& \exp \biggl(-
\lambda\epsilon+ \frac{\lambda^2}{2}\Esp[Z] \biggr)
\\
&=& \exp \biggl(- \frac{(b-a)^2\epsilon^2}{2(n-1)\sigma^2} \biggr) ,
\end{eqnarray*}
where we used the same value $\lambda= \frac{\epsilon}{\Esp[Z]} =
\frac{(b-a)^2\epsilon}{(n-1)\sigma^2}$ 
as in Maurer \cite{Maurer2006}, Theorem~13.

Finally, we have proven that for all $\delta\in[0,1]$, with probability
higher than $1-\delta$,
\begin{eqnarray*}
\max_{1\leq k\leq n}\sigma_{>,k}^2 \leq
\sigma^2 + 2\sqrt{\sigma ^2}\frac{(b-a)(n-1)}{N-n+1}\sqrt{
\frac{\log(1/\delta)}{2(n-1)}} ,%
\end{eqnarray*}
which concludes the proof of \eqref{eqn:sigmaleft}.

We now turn to proving \eqref{eqn:sigmaleft2}.
First, we remark that
\begin{eqnarray*}
\sigma^2_{<,k+1} &\leq& \Esp \bigl[(X_{k+1}-
\mu)^2\vert Z_{k+1},\ldots,Z_{N-1} \bigr]
\\
&=& \Esp \bigl[(X_{k+1}-\mu)^2\vert X_{k+2},
\ldots,X_{N} \bigr]
\\
&=& \Esp \bigl[(Y_{N-k}-\mu)^2\vert Y_1,
\ldots,Y_{N-k-1} \bigr] ,
\end{eqnarray*}
where in the second line we used that $Z_{k+1} = \mu- X_N -\cdots- X_{k+2}$,
and in the third line we used the change of variables
$Y_{u}=X_{N-u+1}$. It follows that
\begin{eqnarray*}
\max_{n\leq k\leq N-1} \sigma^2_{<,k+1} &\leq& \max
_{n\leq
k\leq N-1} \Esp \bigl[(Y_{N-k}-\mu)^2 \vert
Y_1,\ldots,Y_{N-k-1} \bigr]
\\
&=& \max_{1\leq
k\leq N-n}\Esp \bigl[(Y_{k}-\mu)^2
\vert Y_1,\ldots,Y_{k-1} \bigr] .
\end{eqnarray*}
Now $(Y_1,\ldots,Y_{N-n})$ has the same marginal distribution as
$(X_1,\ldots,X_{N-n})$,
so that the proof of \eqref{eqn:sigmaleft} applies and yields the result.
\end{pf*}

We emphasize that we used Hoeffding's reduction Lemma~\ref{lem:reduction} in the proof of Lemma~\ref{lem:sigmaleft}. This
allowed us to apply the key result from Maurer \cite{Maurer2006}. We will
discuss alternatives to this proof in Section~\ref{sec:discussion}. We
can now
state our Bernstein--Serfling bound.
%
\begin{theorem}[(Bernstein--Serfling inequality)]\label{thm:Bern2}
Let $\cX=(x_1,\ldots,x_N)$ be a finite population of $N>1$ real points,
and $(X_1,\ldots,X_n)$ be a list of size $n<N$ sampled without
replacement from $\cX$.
Then, for all $\epsilon>0$ and $\delta\in[0,1]$, the following
concentration inequality holds
%
\begin{eqnarray}\label{eqn:bernSerf1}
\Pr \biggl( \max_{1\leq k\leq n} \frac{\sum_{t=1}^k (X_t - \mu
)}{N-k} \geq
\frac{n\eps}{N - n} \biggr) \leq \exp \biggl[ \frac{-n\epsilon^2/2}{\gamma^2 + \sklfrac
{2}{3}(b-a)\epsilon} \biggr] + \delta,
\end{eqnarray}
where
\[
\gamma^2 = (1-f_{n-1}) \sigma^2 +
f_{n-1}c_{n-1}(\delta) ,
\]
$c_{n}(\delta) = \sigma(b-a)\sqrt{\frac{2\log(1/\delta)}{n}}$, and
$f_{n-1} = \frac{n-1}{N}$.
Similarly, it holds
%
\begin{eqnarray}\label{eqn:bernSerf2}
\Pr \biggl( \max_{n\leq k\leq N-1} \frac{\sum_{t=1}^k (X_t - \mu
)}{k} \geq\eps \biggr)
\leq \exp \biggl[ \frac{- n \epsilon^2/2}{\tilde{\gamma}^2 + \sklfrac
{2}{3}(b-a)\epsilon} \biggr] + \delta,
\end{eqnarray}
where
\[
\tilde{\gamma}^2=(1-f_n) \biggl(\frac{n+1}{n}
\sigma^2 + \frac
{N-n-1}{n} c_{N-n-1}(\delta) \biggr) .
\]
%
\end{theorem}
%
%
\begin{pf}
We first prove \eqref{eqn:bernSerf2}. Applying
Proposition~\ref{prop:Bern} together with Markov's inequality, we
obtain that for all $\lambda,\delta>0$,
%
\begin{eqnarray}\label{eqn:cor1}
\Pr \Biggl( \max_{n\leq k\leq N-1} \frac{\sum_{t=1}^k (X_t - \mu
)}{k} \geq
\frac{\log(1/\delta)}{\lambda} + \lambda \sum_{k=1}^{N-n}
\phi \biggl(\frac{2(b-a) \lambda}{N-k} \biggr)\frac
{\sigma_{<,N-k+1}^2 }{(N -k)^2} \Biggr) \leq
\delta.
\end{eqnarray}
%

Thus, combining equations \eqref{eqn:cor1} and \eqref{eqn:sigmaleft2}
with a
union bound, we get that for all $\lambda>0$ for all $\delta,\delta
'$, with probability higher than $1-\delta-\delta'$, it holds that
\begin{eqnarray*}
&&\max_{n\leq k\leq N-1}\frac{\sum_{t=1}^k (X_t - \mu)}{k}
\\
&&\quad \leq \frac{\log(1/\delta)}{\lambda} + \lambda \sum_{k=1}^{N-n}
\phi \biggl(\frac{2(b-a) \lambda}{N-k} \biggr)\frac
{1}{(N -k)^2} \biggl[
\sigma^2 + \frac{N-n-1}{n+1}c_{N-n-1}\bigl(\delta
'\bigr) \biggr]
\\
&&\quad \leq \frac{\log(1/\delta)}{\lambda} + \frac{\lambda}{n^2} \phi \biggl(\frac{2(b-a) \lambda}{n}
\biggr) \biggl[\sigma^2 + \frac
{N-n-1}{n+1}c_{N-n-1}\bigl(
\delta'\bigr) \biggr] \sum_{k=1}^{N-n}
\frac{n^2}{(N -k)^2}
\\
&&\quad \leq \frac{\log(1/\delta)}{\lambda} + \frac{\lambda}{n^2} \phi \biggl(\frac{2(b-a) \lambda}{n}
\biggr) \biggl[\sigma^2 + \frac
{N-n-1}{n+1}c_{N-n-1}\bigl(
\delta'\bigr) \biggr] (n+1) \biggl(1 - \frac{n}{N} \biggr) ,
\end{eqnarray*}
where we introduced
\[
c_{N-n-1}\bigl(\delta'\bigr) = \sigma(b-a)\sqrt{
\frac{2\log(1/\delta
')}{N-n-1}} ,
\]
where we used in the second line the fact that $\phi$ is nondecreasing
and where we applied \eqref{eqn:sumkn} in the last line.
For convenience, let us now introduce the quantities $f_{n} = \frac
{n}{N}$ and
\begin{eqnarray*}
\tilde{\gamma}^2 &=&(1-f_{n}) \biggl[\sigma^2
+ \frac
{N-n-1}{n+1}c_{N-n-1}\bigl(\delta'\bigr) \biggr] .
\end{eqnarray*}
The previous bound can be rewritten in terms of $\epsilon>0$ and
$\delta'$ only, in the form
%
\begin{eqnarray}
\Pr \biggl( \max_{n\leq k\leq N-1}\frac{\sum_{t=1}^k (X_t - \mu
)}{k} \geq\epsilon \biggr)
\leq\exp \biggl(-\lambda\epsilon+ \frac{\lambda^2(n+1)}{n^2} \phi \biggl(
\frac{2(b-a) \lambda}{n} \biggr)\tilde{\gamma}^2 \biggr) +
\delta' . \label{eqn:invertedBound}
\end{eqnarray}

We now optimize the bound \eqref{eqn:invertedBound} in $\lambda$. Let
us introduce the function
\[
f(\lambda) = -\lambda\epsilon+ \frac{\lambda^2(n+1)}{n^2} \phi \biggl(
\frac{2(b-a) \lambda}{n} \biggr)\tilde{\gamma}^2 ,
\]
corresponding to the term in brackets in \eqref{eqn:invertedBound}. By
definition of $\phi$, it comes
\begin{eqnarray*}
f(\lambda) &=& -\lambda\epsilon+ \frac{\lambda^2}{n^2}\phi \biggl(
\frac{2(b-a) \lambda}{n} \biggr)\tilde{\gamma}^2 (n+1)
\\
&=& -\lambda\epsilon+ \biggl(\exp \biggl(\frac{2(b-a)\lambda}{n} \biggr) - 1 -
\frac{2(b-a)\lambda}{n} \biggr) \frac{\tilde{\gamma
}^2}{4(b-a)^2} (n+1) .
\end{eqnarray*}
Thus, the derivative of $f$ is given by
\begin{eqnarray*}
f'(\lambda) = -\epsilon+ \biggl(\exp \biggl(\frac{2(b-a)\lambda
}{n}
\biggr) -1 \biggr)\frac{\tilde{\gamma}^2 (n+1)}{2(b-a)n} ,
\end{eqnarray*}
and the value $\lambda^\star$ that optimizes $f$ is given by
\begin{eqnarray*}
\lambda^\star= \frac{n}{2(b-a)} \log \biggl(1 + \frac{2(b-a)
\epsilon n}{\tilde{\gamma}^2(n+1)}
\biggr) .
\end{eqnarray*}
Let us now introduce for convenience the quantity $u = \frac{2(b-a)
n}{\tilde{\gamma}^2(n+1)}$.
The corresponding optimal value $f(\lambda^\star)$ is given by
\begin{eqnarray*}
f\bigl(\lambda^\star\bigr) &=& -\epsilon\frac{n}{2(b-a)} \log(1+u
\epsilon) + \frac{\tilde{\gamma}^2}{4(b-a)^2} (n+1) \bigl(u\epsilon- \log (1+u\epsilon) \bigr)
\\
&=& \frac{\tilde{\gamma}^2(n+1)}{4^(b-a)^2} \bigl[ - u\epsilon\log (1+u\epsilon) + u\epsilon-
\log(1+u\epsilon) \bigr]
\\
&=& -\frac{n}{2(b-a) u}\zeta(u\epsilon) ,
\end{eqnarray*}
where we introduced in the last line the function $\zeta(u) =
(1+u)\log(1+u) -u$.
Now, using the identify $\zeta(u) \geq u^2/(2+2u/3)$ for $u\geq0$, we obtain
\begin{eqnarray*}
\Pr \biggl(\max_{n\leq k\leq N-1}\frac{\sum_{t=1}^k (X_t - \mu)}{k} \geq\epsilon \biggr)
&\leq& \exp \biggl( -\frac{n\epsilon}{2(b-a)}\frac{u\epsilon}{2+
2u\epsilon/3} \biggr) +
\delta'
\\
&\leq& \exp \biggl( - \frac{n\epsilon^2}{2\tilde{\gamma}^2(n+1)/n
+ \frac{4}{3}(b-a)\epsilon} \biggr) + \delta' ,
\end{eqnarray*}
which concludes the proof of \eqref{eqn:bernSerf2}. The proof of
\eqref{eqn:bernSerf1} follows the very same lines, simply
using \eqref{eqn:sigmaleft} instead of \eqref{eqn:sigmaleft2}.
\end{pf}

Inverting the bounds of Theorem~\ref{thm:Bern2}, we obtain Corollary~\ref{cor:Bern}.
%
\begin{corollary}\label{cor:Bern}
Let $n\leq N$ and $\delta\in[0,1]$.
With probability larger than $1-2\delta$, it holds that
\begin{eqnarray*}
\frac{\sum_{t=1}^n (X_t - \mu)}{n} \leq \sigma\sqrt{\frac
{2\rho_n\log(1/\delta)}{n}} +
\frac{\kappa_n(b-a)\log(1/\delta)}{n},
\end{eqnarray*}
where we remind the definition of $\rho_n$ \eqref{eqn:rhon}
\[
\rho_n = %
\cases{ (1-f_{n-1}) & \quad \mbox{if} $n\leq
N/2$,
\cr
(1-f_{n}) (1+1/n) & \quad \mbox{if} $n > N/2$, } %
\]
and where we introduced the quantity
%
\begin{eqnarray}\label{eqn:kappan}
\kappa_n = %
\cases{ \displaystyle \frac{4}{3} + \sqrt{
\displaystyle \frac{f_n}{g_{n-1}}} & \quad \mbox{if} $n\leq N/2$,
\cr
\displaystyle \frac{4}{3}+
\sqrt{g_{n+1}(1-f_n)}& \quad \mbox{if} $n > N/2$, }
\end{eqnarray}
with $f_n = n/N$ and $g_n = N/n-1$.
\end{corollary}
%
%
\begin{pf}
Let $\delta,\delta'\in[0,1]$. From \eqref{eqn:bernSerf1} in
Theorem~\ref{thm:Bern2}, it comes that, with probability higher than
$1-\delta-\delta'$,
\begin{eqnarray*}
\frac{\sum_{t=1}^n (X_t - \mu)}{N-n} \leq\epsilon_\delta,\qquad  \mbox {where }
\gamma^2 + B \frac{N-n}{n}\epsilon_\delta=
\frac
{(N-n)^2}{2n\log(1/\delta)}\epsilon_\delta^2 ,
\end{eqnarray*}
where we introduced for convenience $B = \frac{2}{3}(b-a)$ and
\begin{eqnarray*}
\gamma^2 = (1-f_{n-1}) \sigma^2 +
f_{n-1}\sigma(b-a)\sqrt{\frac
{2\log(1/\delta')}{n-1}} .
\end{eqnarray*}
Solving this equation in $\epsilon$ leads to
\begin{eqnarray*}
\epsilon_\delta&=& n\log(1/\delta)\frac{B\sklvfrac{N-n}{n} + \sqrt
{B^2 (\vfrac{N-n}{n} )^2 + 4 \sklafrac{(N-n)^2}{2n\log(1/\delta
)}\gamma^2 }}{(N-n)^2}
\\
&=& \frac{1}{N-n} \bigl(\sqrt{B^2 \log(1/\delta)^2
+ 2\gamma^2\log (1/\delta)n} + B\log(1/\delta) \bigr)
\\
&\leq& \frac{n}{N-n} \biggl( \sqrt{\frac{2\gamma^2\log(1/\delta
)}{n}} +
\frac{2B\log(1/\delta)}{n} \biggr) .
\end{eqnarray*}

On the other hand, following the same lines but starting from \eqref
{eqn:bernSerf2} in Theorem~\ref{thm:Bern2}, it holds that, with
probability higher than $1-\delta-\delta'$,
\begin{eqnarray*}
\frac{\sum_{t=1}^n (X_t - \mu)}{n} \leq\sqrt{\frac{2\tilde
\gamma^2\log(1/\delta)}{n}} + \frac{2B\log(1/\delta)}{n} ,
\end{eqnarray*}
where we introduced this time
\begin{eqnarray*}
\tilde\gamma^2 = (1-f_{n}) \biggl((1+1/n)
\sigma^2 + \frac
{N-n-1}{n}\sigma(b-a)\sqrt{\frac{2\log(1/\delta')}{N-n-1}}
\biggr) .
\end{eqnarray*}
Finally, we note that
\begin{eqnarray*}
\sqrt{\tilde\gamma^2} &\leq&
\sqrt{(1-f_{n}) (1+1/n)} \biggl(\sigma+ \frac{N-n-1}{n+1}(b-a)\sqrt{
\frac{\log(1/\delta')}{2(N-n-1)}} \biggr) .
\end{eqnarray*}

Thus, when $n\leq N/2$, we deduce that for all $1\leq n\leq N-1$, with
probability higher than $1-2\delta$, it holds
\begin{eqnarray*}
\frac{\sum_{t=1}^n (X_t - \mu)}{n} &\leq& \sqrt{1-f_{n-1}} \biggl(\sigma\sqrt{
\frac{2\log(1/\delta)}{n}} + \frac{n-1}{N-n+1}\frac
{(b-a)\log(1/\delta)}{\sqrt{n(n-1)}} \biggr)
\\
&&{} + \frac{2B\log(1/\delta)}{n}
\\
&\leq& \sigma\sqrt{\frac{2(1-f_{n-1})\log(1/\delta)}{n}} + \frac
{(b-a)\log(1/\delta)}{n} \biggl(
\frac{4}{3}+\sqrt{\frac
{n(n-1)}{N(N-n+1)}} \biggr) ;
\end{eqnarray*}
whereas when $N>n>N/2$, it holds, with probability higher than
$1-2\delta$, that
\begin{eqnarray*}
\frac{\sum_{t=1}^n (X_t - \mu)}{n} &\leq& \sqrt {(1-f_{n}) (1+1/n)} \biggl(\sigma
\sqrt{\frac{2\log(1/\delta)}{n}} + \frac{N-n-1}{n+1}\frac{(b-a)\log(1/\delta)}{\sqrt{n(N-n-1)}} \biggr)
\\
&&{}+ \frac{2B\log(1/\delta)}{n}
\\
&\leq& \sigma\sqrt{\frac{2(1-f_{n})(1+1/n)\log(1/\delta)}{n}}
\\
&&{}+ \frac{(b-a)\log(1/\delta)}{n} \biggl(\frac{4}{3}+ \sqrt{\frac
{(N-n-1)(N-n)}{(n+1)N}}
\biggr) .
\end{eqnarray*}
Finally we note that when $n=N$, $g_{n+1}(1-f_n)=0$ and $\rho_n = 0$.
So the bound is still\linebreak[4]  satisfied.
\end{pf}

\section{An empirical Bernstein--Serfling inequality}\label{sec:empBern}

In this section, we derive a practical version of Theorem~\ref{thm:Bern2}
where the variance $\sigma^2$ is replaced by an estimate.
A natural (biased) estimator is given by
%
\begin{eqnarray}\label{eqn:defEmpVar}
\hat\sigma_n^2 = \frac{1}{n}\sum
_{i=1}^n(X_i - \hat
\mu_n)^2 = \frac{1}{n^2}\sum
_{i,j=1}^n \frac{(X_i-X_j)^2}{2}, \qquad \mbox{where } \hat
\mu_n = \frac{1}{n}\sum_{i=1}^n
X_i .
\end{eqnarray}
We also define, for notational convenience,
the quantity $\hat\sigma_n = \sqrt{\hat\sigma_n^2}$.

Before proving our empirical Bernstein--Serfling inequality, we first
need to control the error between $\hat\sigma_n$ and
$\sigma$. For instance, in the standard case of sampling with
replacement, it can be shown (Maurer and Pontil \cite{MaurerP09}) that, for all $\delta
\in[0,1]$,
\begin{eqnarray*}
\Pr \biggl(\sigma\geq\frac{n}{n-1}\hat\sigma_n + (b-a)\sqrt{
\frac
{2\ln(1/\delta)}{n-1}} \biggr) \leq\delta.
\end{eqnarray*}
We now show an equivalent result in the case of sampling without replacement.
%
\begin{lemma}\label{lem:sigmaright}
When sampling without replacement from a finite population
$\cX=(x_1,\ldots,x_N)$ of size $N$, with range $[a,b]$ and variance
$\sigma^2$,
the empirical variance $\hat\sigma_n^2$ defined in \eqref
{eqn:defEmpVar} using $n<N$ samples satisfies the following
concentration inequality (using the notation of Corollary~\ref{cor:Hoeff})
\begin{eqnarray*}
\Pr \biggl(\sigma\geq\hat\sigma_n + (b-a) (1+\sqrt{1+\rho
_n} )\sqrt{\frac{\log(3/\delta)}{2n}} \biggr) \leq\delta.
\end{eqnarray*}
\end{lemma}
%
\begin{remark}
We conjecture that it is possible, at the price of a more
complicated analysis, to reduce the term
$(1+\sqrt{1+\rho_n})$
to $\sqrt{4\rho_n}$, which would then be consistent with the
analogous result for sampling with replacement in Maurer and Pontil \cite{MaurerP09}.
We further
discuss this technically involved improvement in Section~\ref{sec:discussion}.
\end{remark}
\begin{pf*}{Proof of Lemma~\ref{lem:sigmaright}}
In order to prove Lemma~\ref{lem:sigmaright}, we again use Lemma~\ref
{lem:reduction},
which allows us to relate the concentration of the quantity $V_n =
\frac{1}{n}\sum_{i=1}^n (X_i -\mu)^2$ to that of
its equivalent
\[
\tilde V_n=\tilde V_n(\mathbf{Y}_{1:n}) =
\frac{1}{n}\sum_{i=1}^n
(Y_i -\mu)^2 ,
\]
where the $Y_i$s are drawn from $\cX$ \textit{with replacement}.
Let us introduce the notation $Z = \frac{n}{(b-a)^2}\tilde V_n(\mathbf
{Y}_{1:n})$.
We know from the proof of Lemma~\ref{lem:sigmaleft} that $Z$ satisfies
the conditions of application of Maurer \cite{Maurer2006}, Theorem~13.
Let us also introduce for convenience the constant $\lambda= -\frac
{\epsilon}{\Esp[Z]} = -\frac{(b-a)^2\epsilon}{n\sigma^2}$.
Using these notations, it comes\vspace*{-2pt}
\begin{eqnarray*}
\Pr \biggl( \sigma^2 - V_n \geq\frac{(b-a)^2}{n}
\epsilon \biggr) &\leq& \Esp \biggl[\exp \biggl(-\lambda \biggl(
\frac{n}{(b-a)^2}\sigma^2 - \frac{n}{(b-a)^2}V_n -
\epsilon \biggr) \biggr) \biggr]
\\[-1pt]
&\leq& \Esp \bigl[\exp \bigl(-\lambda \bigl(\Esp[Z] - Z -\epsilon \bigr) \bigr)
\bigr]
\\[-1pt]
&\leq& \exp \biggl(\lambda\epsilon+ \frac{\lambda^2}{2}\Esp[Z] \biggr)
\\[-1pt]
&=& \exp \biggl( -\frac{(b-a)^2\epsilon^2}{2n\sigma^2} \biggr) .
\end{eqnarray*}
The first line results of the application of Markov's inequality.
The second line follows from the application of Lemma~\ref
{lem:reduction} to $X'_i = (X_i-\mu)^2$ and
$f(x) = \exp (-\lambda\frac{n}{(b-a)^2}x )$.
The last steps are the same as in the proof of Lemma~\ref{lem:sigmaleft}.

So far, we have shown that, with probability at least $1-\delta$,\vspace*{-2pt}
%
\begin{eqnarray}\label{eqn:sigVn}
\sigma^2- 2\sqrt{\sigma^2}(b-a)\sqrt{
\frac{\log(1/\delta)}{2n }} \leq V_n .
\end{eqnarray}
Let us remark that\vspace*{-2pt}
\begin{eqnarray*}
\frac{1}{n}\sum_{i=1}^n
(X_i -\mu)^2- \frac{1}{n}\sum
_{i=1}^n (X_i -\hat
\mu_n)^2 &=& (\hat\mu_n-\mu)^2 ,
\end{eqnarray*}
that is, $V_n = (\hat\mu_n-\mu)^2 + \hat\sigma_n^2$. In order to
complete the proof, we thus resort twice to Theorem~\ref{thm:Hoeff} to
obtain that, with probability higher than $1-\delta$,
it holds\vspace*{-2pt}
%
\begin{eqnarray}\label{eqn:hoefsqr}
(\hat\mu_n-\mu)^2 \leq(b-a)^2
\frac{\rho_n\log(2/\delta)}{2n} .
\end{eqnarray}
Combining equations \eqref{eqn:sigVn} and \eqref{eqn:hoefsqr} with a
union bound argument yields that, with
probability at least $1-\delta$,\vspace*{-2pt}
\begin{eqnarray*}
\hat\sigma_n^2 &\geq& \sigma^2- 2\sqrt{
\sigma^2}\sqrt {(b-a)^2\frac{\log(3/\delta)}{2n }} -
(b-a)^2\frac{\rho_n\log
(3/\delta)}{2n}
\\[-1pt]
&=& \biggl(\sigma- \sqrt{(b-a)^2\frac{\log(3/\delta)}{2n }}
\biggr)^2 - (b-a)^2 (1+\rho_n )
\frac{\log(3/\delta)}{2n} .
\end{eqnarray*}
Finally, we obtain\vspace*{-1pt}
\begin{eqnarray*}
\Pr \biggl(\sigma\geq\hat\sigma_n + (1+\sqrt{1+\rho_n} )
\sqrt{(b-a)^2\frac{\log(3/\delta)}{2n}} \biggr) \leq\delta.
\end{eqnarray*}
\upqed
\end{pf*}

Eventually, combining Theorem~\ref{thm:Bern2} and
Lemma~\ref{lem:sigmaright} with a union bound argument, we finally
deduce the following result.
%
\begin{theorem}[(An empirical Bernstein--Serfling inequality)]\label{thm:EmpBern}
Let $\cX=(x_1,\ldots,x_N)$ be a finite population of $N>1$ real points,
and $(X_1,\ldots,X_n)$ be a list of size $n\leq N$ sampled without
replacement from $\cX$.
Then for all $\delta\in[0,1]$, with probability larger than
$1-5\delta$, it holds
\[
\frac{\sum_{t=1}^n (X_t - \mu)}{n} \leq\hat \sigma_n\sqrt{\frac{2\rho_n\log(1/\delta)}{n}} +
\frac{\kappa(b-a)\log(1/\delta)}{n} ,
\]
where we remind the definition of $\rho_n$ \eqref{eqn:rhon}
\[
\rho_n = %
\cases{ \biggl(1-\displaystyle \frac{n-1}{N}\biggr) &
\quad \mbox{if} $n\leq N/2$,
\cr
\biggl(1-\displaystyle \frac{n}{N}\biggr) (1+1/n) &\quad  \mbox{if} $n
> N/2$ , } %
\]
and $\kappa= \frac{7}{3}+\frac{3}{\sqrt{2}}$.
\end{theorem}
%
%
\begin{remark}
First, Theorem~\ref{thm:EmpBern} has the familiar form of
Bernstein bounds. The alternative definition of $\rho_n$ guarantees
that we get
the best reduction out of the no replacement setting. In particular,
when $n$ is large, the factor $(1-f_n)$ replaces $(1-f_{n-1})$ and the
corresponding factor
eventually equals $0$ when $n=N$, a feature that was missing in
Proposition~\ref{prop:Hoeff}. Second, the constant
$\kappa$ is to relate to the constant $7/3$ in
Maurer and Pontil \cite{MaurerP09}, Theorem~11, for sampling
with replacement.
\end{remark}
\begin{pf*}{Proof of Theorem~\ref{thm:EmpBern}}
First, by application of Corollary~\ref{cor:Bern}, it holds for all
$\delta\in[0,1]$ that, with probability
higher than $1-2\delta$,
\begin{eqnarray*}
\frac{\sum_{t=1}^n (X_t - \mu)}{n} &\leq& \sigma\sqrt{\frac
{2\rho_n\log(1/\delta)}{n}} +
\frac{\kappa_n(b-a)\log(1/\delta)}{n},
\end{eqnarray*}
where we remind the definition of $\rho_n$ \eqref{eqn:rhon}
\[
\rho_n = %
\cases{ (1-f_{n-1}) & \quad \mbox{if} $n\leq
N/2$,
\cr
(1-f_{n}) (1+1/n) & \quad \mbox{if} $n > N/2$, } %
\]
and the definition of $\kappa_n$ \eqref{eqn:kappan}
\[
\kappa_n = %
\cases{ \displaystyle \frac{4}{3} + \sqrt{
\displaystyle \frac{f_n}{g_{n-1}}} & \quad \mbox{if} $n\leq N/2$,
\cr
\displaystyle \frac{4}{3}+
\sqrt{g_{n+1}(1-f_n)}& \quad \mbox{if} $n > N/2$. } %
\]
%
We then apply Lemma~\ref{lem:sigmaright} to get that, with probability
higher than $1-5\delta$,
if $n\leq N/2$, then
%
\begin{eqnarray}\label{eqn:empBernTool1}
\frac{\sum_{t=1}^n (X_t - \mu)}{n} &\leq& \sqrt{\hat\sigma _n^2}
\sqrt{\frac{2\log(1/\delta)}{n}} \sqrt{1-f_{n-1}}
\nonumber
\\
&&{} + \frac{(b-a)\log(1/\delta)}{n} \biggl(\frac{4}{3}+\sqrt{\frac{f_n}{g_{n-1}}}
\\
&&\hphantom{{} + \frac{(b-a)\log(1/\delta)}{n} \biggl(}{}+(1+ \sqrt{2-f_{n-1}})\sqrt{1-f_{n-1}} \biggr) ,\nonumber\end{eqnarray}
and if $n>N/2$, then
%
\begin{eqnarray}\label{eqn:empBernTool2}
&&\frac{\sum_{t=1}^n (X_t - \mu)}{n} \nonumber\\
&&\quad \leq \sqrt{\hat\sigma _n^2}
\sqrt{\frac{2\log(1/\delta)}{n}}\sqrt {(1-f_{n}) (1+1/n)}
\nonumber
\\[-8pt]\\[-8pt]
&&\qquad {} +\frac{(b-a)\log(1/\delta)}{n} \biggl(\frac{4}{3}+ \sqrt {g_{n+1}(1-f_n)}\nonumber
\\
&&\hphantom{\qquad {} +\frac{(b-a)\log(1/\delta)}{n} \biggl(}{}+ \sqrt{(1-f_{n}) (1+1/n)} \bigl(1+ \sqrt{1+ (1-f_{n})
(1+1/n)} \bigr) \biggr) . \nonumber
\end{eqnarray}

We now simplify this result.
%
Assume first that $n\leq N/2$. We thus get
\[
\frac{f_n}{g_{n-1}} \leq\frac{1}{2g_{n-1}} = \frac{n-1}{2(N-n+1)} \leq
\frac{1}{2},
\]
so that we deduce
%
\begin{eqnarray}\label{eqn:empBernTerm1}
 \frac{4}{3}+(1+ \sqrt{2-f_{n-1}})\sqrt{1-f_{n-1}}+
\sqrt{\frac
{f_n}{g_{n-1}}} \leq2+ \frac{1}{3}+\sqrt{2}+\frac{1}{\sqrt{2}}
.
\end{eqnarray}
Assume now that $n> N/2$. In this case, it holds
\begin{eqnarray*}
g_{n+1}(1-f_n)& =& \frac{N-n-1}{n+1}\frac{N-n}{N}
\leq\frac{N-n}{N} \leq \frac{1}{2},
\\
(1-f_{n}) (1+1/n)& = &\biggl(1-\frac{n}{N} \biggr) (1+1/n) \leq
\frac
{1}{2} \biggl(1+\frac{2}{N} \biggr) ,
\end{eqnarray*}
so that we deduce, since $N\geq{2}$,
\begin{eqnarray}\label{eqn:empBernTerm2}
\frac{4}{3}+\sqrt{g_{n+1}(1-f_n)} +
\sqrt{(1-f_{n}) (1+1/n)} (1+ \sqrt{2-f_{n-1}}) \leq2+
\frac{1}{3} + \frac{1}{\sqrt{2}} +\sqrt{2} .
\end{eqnarray}
%

Respectively combining \eqref{eqn:empBernTerm1} and \eqref
{eqn:empBernTerm2} with
equations \eqref{eqn:empBernTool1} and \eqref{eqn:empBernTool2}
concludes the proof.
\end{pf*}

\section{Discussion}
\label{sec:discussion}
In this section, we discuss the bounds of Theorem~\ref{thm:Bern2}
and Theorem~\ref{thm:EmpBern} from the perspective of both theory and
application.
%
\begin{figure}

\includegraphics{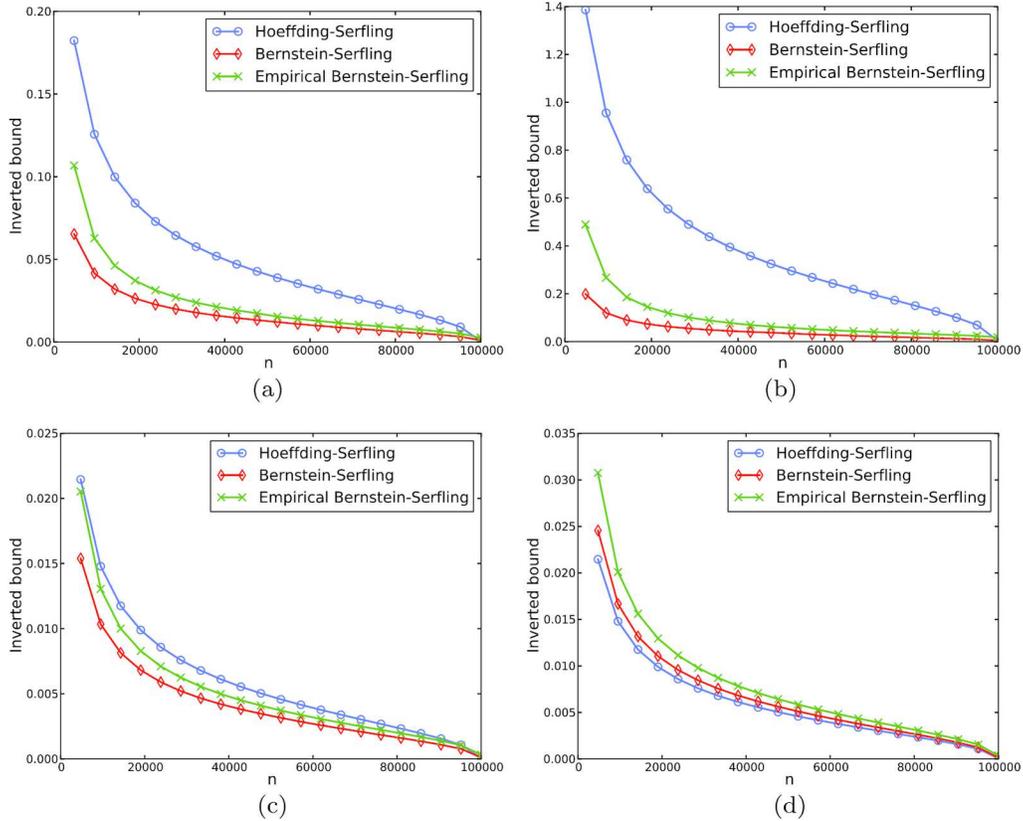}

\caption{Comparing the bounds of Corollaries \protect\ref{cor:Hoeff} and
\protect\ref{cor:Bern}, and Theorem \protect\ref{thm:EmpBern}. $\cX$ is here a
sample from each of the four distributions written below each plot,
of size $N=10^6$. Unlike Figure~\protect\ref{fig:plots}, as $n$ increases,
we keep sampling here without replacement until exhaustion.
(a) Gaussian $\cN(0,1)$. (b) Log-normal $\ln\cN(1,1)$. (c) Bernoulli
$\cB(0.1)$.
(d) Bernoulli $\cB(0.5)$.}\label{fig:empBern}
\end{figure}

First, both bounds involve either the factor $1-f_{n-1}$ or $1-f_n$,
thus leading to a dramatic improvement on the usual Bernstein or
empirical Bernstein bounds, which do not make use of the no
replacement setting. This is crucial, for
instance, when the user needs to rapidly compute an empirical mean
from a large number of samples up to some precision level. To better
understand the improvement of Serfling bounds, we plot in
Figure~\ref{fig:empBern} the bounds of
Corollaries \ref{cor:Hoeff} and \ref{cor:Bern}, and Theorem~\ref
{thm:EmpBern} for an example where $\cX$ is
a sample of size $N=10^6$ from each of the following four
distributions: unit centered Gaussian, log-normal with parameters
$(1,1)$, and Bernoulli with parameter 1$/$10 and 1$/$2. As $n$ increases, we
keep sampling without replacement from $\cX$ until
exhaustion, and report the corresponding bounds. Note that all our
bounds have their leading term exactly equal to zero when $n=N$, though our
Hoeffding--Serfling bound only is exactly zero. In all experiments, the
loss of tightness as a result of using the empirical variance is
small. Our empirical Bernstein--Serfling demonstrates here a dramatic
improvement on the Hoeffding--Serfling bound of
Corollary~\ref{cor:Hoeff} in Figures~\ref{fig:empBern}(a) and
\ref{fig:empBern}(b). A slight improvement is demonstrated in
Figure~\ref{fig:empBern}(c) where the standard deviation of
$\cX$ is roughly a third of the range. Finally, Bernstein--Serfling itself
does not improve on Hoeffding--Serfling in
Figure~\ref{fig:empBern}(d), where the standard deviation is
roughly half of the range, again indicating that Bernstein bounds are
not uniformly better than Hoeffding bounds.

There is a number of nontrivial applications of our
bounds. \emph{Scratch games}, for instance, were introduced in
F{\'e}raud and Urvoy \cite{feraud2013exploration} as
a variant of the
multi-armed bandit problem, to model two real world problems: selecting
ads to display on web pages and optimizing e-mailing campaigns.
In particular, F{\'e}raud and Urvoy \cite{feraud2013exploration} discuss practical
situations where an upper confidence bound algorithm based on a
Hoeffding--Serfling inequality
outperforms a standard algorithm based on Hoeffding's
inequality. Similar improvements should appear in practice when using
our empirical Bernstein--Serfling inequality.
As another application, our results could be useful in
optimization. The stochastic dual-coordinate ascent algorithm (SDCA;
Shalev-Shwartz and Zhang \cite{ShalevShwartz013}) is a
state-of-the-art optimization algorithm
used in machine learning. Shalev-Shwartz and Zhang \cite{ShalevShwartz013} introduce a variant
of SDCA called SDCA-Perm, which -- unlike SDCA -- relies on sampling
without replacement,
and achieves better empirical performance than SDCA. However, the
analysis in Shalev-Shwartz and Zhang \cite{ShalevShwartz013}
does not cover SDCA-Perm. We
believe that the use of Serfling bounds is an appropriate tool for
that purpose.

To conclude, we discuss potential improvements of our bounds. A careful
look at Lemmas \ref{lem:sigmaleft}
and \ref{lem:sigmaright} indicates that our bounds may be further
improved, though at the price of a more intricate
analysis. Indeed, these two lemmas both resort to Hoeffding's
reduction Lemma~\ref{lem:reduction}, in order to be able to apply
concentration results known for self-bounded random variables to the
setting of sampling without replacement. As a result, we lose here a
potential factor $\rho_n$ for the confidence bound around the
variance, and we conjecture that the term $1+\sqrt{1+\rho_n}$ in
Lemma~\ref{lem:sigmaright} could ultimately be replaced with
$2\sqrt{\rho_n}$. A natural tool for this
would be a dedicated
\textit{tensorization inequality} for the entropy in the case of
sampling without replacement\vadjust{\goodbreak}
(Boucheron, Lugosi and Massart \cite{BoLuMa13}, Maurer
\cite{Maurer2006}, Bousquet \cite{Bousquet03}).
Indeed, it is not difficult to
show that $\hat\sigma_{n}^2$ satisfies a self-bounded property
similar to
that of Maurer and Pontil \cite{MaurerP09}, Theorem~11,
involving the factor $\rho_n$.
Thus, in order to be able to get a version of Maurer and Pontil
\cite{MaurerP09}, Theorem~11, in our setting, a specific so-called tensorization
inequality would be enough. Unfortunately, we are unaware of the
existence of such an inequality for sampling without replacement,
where the samples are strongly dependent. We are also unaware of any
tensorization inequality designed for $U$-statistics, which could be
another possible way to get the desired result. Although we believe
this is possible, developing such tools goes beyond the scope of this
paper, and the current results of Theorem~\ref{thm:Bern2} and Theorem~\ref{thm:EmpBern} are already appealing without resorting to further
technicalities,
which would only affect second-order terms in the end.


\section*{Acknowledgements}
This work was supported by both the \textit{2020 Science} program,
funded by EPSRC grant number EP/I017909/1, and the Technion.



\printhistory
\end{document}